\documentclass[12pt]{article}
\usepackage{epsf,epsfig,amsmath,amsfonts}
\usepackage{graphicx}

\newtheorem{thm}{Theorem}
\newtheorem{prop}[thm]{Proposition}
\newtheorem{ex}{Example}

\newtheorem{lem}[thm]{Lemma}
\newtheorem{de}{Definition}
\newtheorem{rem}{Remark}
\newtheorem{ass}{Assumption}

\newcommand{\be}{\begin{eqnarray*}}
\newcommand{\ben}{\begin{eqnarray}}
\newcommand{\ee}{\end{eqnarray*}}
\newcommand{\een}{\end{eqnarray}}
\newcommand{\aln}{\begin{align}}
\newcommand{\ealn}{\end{align}}

\newcommand{\half}{{\textstyle \frac{1}{2}}}

\newcommand{\drq}{{\mathcal D}_r Q}
\newcommand{\drg}{{\mathcal D}_r G}

\newcommand{\F}{{\cal F}}
\newcommand{\e}{{\mathbb E}}
\newcommand{\p}{{\mathbb P}}
\def\Om{\Omega}
\def\R{\mathbb R}

\begin{document}

\title{A Singular Control Model with Application to the
Goodwill Problem\footnote{Research supported by EPSRC grant
no.\,GR/S22998/01}}

\author{
{\sc Andrew Jack}\footnote{Department of Mathematics,
King's College London, The Strand, London WC2R 2LS, UK,
\texttt{andrew.j.f.jack@googlemail.com}}, \
{\sc Timothy C.\,Johnson}\footnote{Department of Actuarial
Mathematics and Statistics, School of Mathematical and
Computer Sciences, Heriot-Watt University, Edinburgh
EH14 4AS, UK, \texttt{t.c.johnson@hw.ac.uk}} \ and
{\sc Mihail Zervos}\footnote{Department of Mathematics, London
School of Economics, Houghton Street, London WC2A 2AE, UK,
\texttt{m.zervos@lse.ac.uk}}}

\maketitle

\begin{abstract}
We consider a stochastic system whose uncontrolled state dynamics are
modelled by a general one-dimensional It\^{o} diffusion.
The control effort that can be applied to this system takes the form
that is associated with the so-called monotone follower problem of
singular stochastic control.
The control problem that we address aims at maximising a performance
criterion that rewards high values of the utility derived from the
system's controlled state but penalises any expenditure of control
effort.
This problem has been motivated by applications such as the so-called
goodwill problem in which the system's state is used to represent the
image that a product has in a market, while control expenditure is
associated with raising the product's image, e.g., through
advertising.
We obtain the solution to the optimisation problem that we consider
in a closed analytic form under rather general assumptions.
Also, our analysis establishes a number of results that are concerned
with analytic as well as probabilistic expressions for the first
derivative of the solution to a second order linear non-homogeneous
ordinary differential equation.
These results have independent interest and can potentially be of use
to the solution of other one-dimensional stochastic control problems.
\end{abstract}

%==================================================================
\section{Introduction}

We consider a stochastic system whose state is modelled by the
controlled, one-dimensional, positive It\^{o} diffusion
\be
dX_t = b(X_t) \, dt + dZ_t + \sigma (X_t) \, dW_t , \quad X_0 = x
> 0 ,
\ee
where $W$ is a standard one-dimensional Brownian motion and
the controlled process $Z$ is a c\`{a}gl\`{a}d increasing
process.
The objective of the optimisation problem is to maximise the
performance criterion
\ben
J_x(Z) = \limsup _{T \rightarrow \infty} \e \left[ \int _0^T
e^{-\Lambda _t^{r (X)}} h(X_t) \, dt - \int _0^T e^{-\Lambda
_t^{r (X)}} k(X_t) \circ dZ_t \right] , \label{J-intro}
\een
over all admissible choices of $Z$, where
\begin{gather}
\Lambda _t^{r (X)} = \int _0^t r(X_\upsilon) \, d\upsilon ,
\nonumber \\
\intertext{and}
\int _0^T e^{-\Lambda _t^{r (X)}} k(X_t) \circ dZ_t
= \int_0^T e^{-\Lambda _t^{r (X)}} k(X_t) \, dZ_t^c +
\sum _{0 \leq t \leq T} \int _0^{\Delta Z_{t}} e^{-\Lambda _t
^{r(X)}} k(X_t+s) \, ds . \label{Z-int}
\end{gather}
This stochastic control problem has been motivated by the following
application.
Consider a company marketing a given product.
This product has an image in a given market where it is being sold
that evolves randomly over time.
We use the random variable $X_t$ to model the product's image at time
$t$, for $t \geq 0$.
The company marketing the product can raise its image by means of
costly interventions such as advertising.
We model the effect of these actions by means of the controlled
process $Z$.
The company's objective is to maximise the expected discounted utility
derived from the product's image minus the expected discounted
``dis-utility'' resulting from intervention costs, which is reflected
in the structure of the performance criterion given by (\ref{J-intro}).

Optimal control problems motivated by applications such as the one
discussed briefly above have a long history and can be traced back to
Nerlove and Arrow~\cite{NA}, who use deterministic dynamics to model
the evolution of the product's image and consider a multi-objective
performance criterion.
Since then, deterministic optimal control models addressing this type
of applications have attracted significant interest (see Buratto and
Viscolani~\cite{BV} and references therein).
Several models in which the product's image evolves randomly over
time, which are more realistic, and result in stochastic optimisation
problems have also been studied in the literature (see
Marinelli~\cite{Mar} and references therein).

Singular stochastic control was introduced by Bather and
Chernoff~\cite{BC} who considered a simplified model of spaceship
control.
In their seminal paper, Bene\v{s}, Shepp and Witsenhausen~\cite{BSW}
were the first to solve rigorously an example of a finite-fuel singular
control problem.
Since then, the area has attracted considerable interest in the
literature.
Alvarez \cite{A1, A2},
Chow, Menaldi and Robin \cite{CMR},
Davis and Zervos \cite{DZ},
Fleming and Soner \cite[Chapter~VIII]{FS},
Harrison and Taksar \cite{HT},
Jacka \cite{J1, J2},
Karatzas \cite{Ka},
K{\o}bila \cite{Ko},
Ma \cite{Ma},
{\O}ksendal \cite{O},
Shreve, Lehoczky and Gavers \cite{SLG},
Soner and Shreve \cite{SS},
Sun \cite{Su},
Zhu \cite{Z},
provide an incomplete list, in alphabetical order, of further
important contributions.
Other related contributions include
Menaldi and Robin~\cite{MR},
Weerasinghe~\cite{W},
and Jack and Zervos~\cite{JZ} who solve singular control problems with
long-term average rather than discounted criteria.
With regard to the structure of the performance criterion that we
consider, penalising the expenditure of control effort by means of
integrals as in (\ref{Z-int}) was introduced by Zhu~\cite{Z} and was
later adopted by Davis and Zervos~\cite{DZ} and Jack and
Zervos~\cite{JZ}.

We solve the problem that we consider by constructing a solution
to the associated Hamilton-Jacobi-Bellman (HJB) equation in closed
analytic form, under general assumptions.
This is possible in the generality that we consider because the control
problem's state space is one-dimensional.
Explicitly solvable control problems have attracted significant
interest in the literature for several reasons.
First, some of them, such as the one that we solve here, are motivated
by real-life applications.
Second, they reveal the qualitative nature of the associated optimal
control tactics, and they provide special cases that can be used to
assess the efficiency of numerical techniques devised to address more
complex problems, which is a major issue.
The majority of such control problems assume that the system's
uncontrolled dynamics are modelled by a Brownian motion with drift or
a geometric Brownian motion.
To the best of our knowledge, Alvarez~\cite{A2} and Jack and
Zervos~\cite{JZ} are the only references in the singular stochastic
control literature in which closed-form solutions are derived when the
system's uncontrolled dynamics are modelled by a general
one-dimensional It\^{o} diffusion.
The latter reference considers a long-term average rather than an
expected discounted performance criterion such as the one given by
(\ref{J-intro}).
On the other hand, Alvarez~\cite{A2} assumes that $h \equiv 0$ and
that $k<0$ is constant.
The introduction of a non-trivial running payoff function $h$ and a
non-trivial running control expenditure function $k$ of sign opposite
to the one of $h$ gives rise to a genuinely new problem that involves
new analysis.
At this point, we should mention that we are not aware of any
reference addressing the problem that we solve, even for simple
stochastic dynamics.

The paper is organised as follows.
Section~\ref{pr-form} is concerned with the formulation of the
singular stochastic control problem that we solve.
In this section, we also develop all of the assumptions on the problem
data that we make in the paper, and we prove that these are sufficient
for our optimisation problem to be well-posed.
Section~\ref{sec:ODE} is concerned with properties of the solution to
a non-homogeneous second-order linear ordinary differential equation
that plays an important role in our analysis.
All of the claims that we make there without proof are standard and
can be found in several references, including Feller~\cite{F},
Breiman~\cite{B}, Mandl~\cite{Man}, It\^{o} and McKean~\cite{IM},
Karlin and Taylor~\cite{KT}, Rogers and Williams~\cite{RW}, and
Borodin and Salminen~\cite{BS}.
The results that we prove are new and have independent interest
because they can be of use in the solution of other stochastic control
problems.
We should note that these results can be generalised with little effort
to account for the case where the underlying diffusion can take values
in an arbitrary interval $\mathcal{I} \in \mathbb{R}$ rather than in
$]0,\infty[$.
Also, some of the assumptions made here can be relaxed.
However, we decided against such relaxations because they would
complicate significantly the exposition.
In Section~\ref{sec:soltn}, we solve the stochastic control problem
that we consider.
Finally, Section~\ref{sec:examples} is concerned with special cases
that arise when $h$ is a power utility function, $k$ and $r$ are
constants, and the uncontrolled state space dynamics are modelled by a
geometric Brownian motion (Section~\ref{sec:GBMex}) or a
mean-reverting square-root process such as the one in the
Cox-Ingersoll-Ross interest rate model (Section~\ref{sec:CIRex}).

%===============================================================
\section{The singular stochastic control problem}
\label{pr-form}

We fix a filtered probability space $(\Om , \F , \F_t ,P)$ satisfying
the usual conditions and carrying a standard one-dimensional
$(\F_t)$-Brownian motion $W$.
We consider a stochastic system whose uncontrolled dynamics
are modelled by the It\^{o} diffusion associated with the stochastic
differential equation (SDE)
\ben
dX_t^0 = b(X_t^0) \, dt + \sigma (X_t^0) \, dW_t , \quad X_0^0 =
x > 0 , \label{SDE-X}
\een
and we make the following assumption.

\begin{ass} {\rm
The functions $b , \sigma : \, ]0,\infty [ \, \rightarrow \R$ are
$C^1$, $\sigma '$ is locally Lipschitz, and $\sigma ^2 (x) > 0$, for
all $x>0$.
} \mbox{}\hfill$\Box$ \label{A1} \end{ass}
This assumption implies that (\ref{SDE-X}) has a unique strong
solution.
It also implies that the scale function $p_{X^0}$ and the
speed measure $m_{X^0}$ given by
\begin{gather}
p_{X^0} (c) = 0 , \quad p_{X^0} '(x) = \exp \left( - 2 \int _c^x
\frac{b(s)} {\sigma ^2(s)} \, ds \right) , \label{scale} \\
\intertext{and}
m_{X^0} (dx) = \frac{2}{\sigma ^2 (x) p_{X^0} ' (x)} \, dx ,
\label{speed}
\end{gather}
respectively, for some $c>0$ fixed, are well-defined.
Additionally, we assume that the solution to (\ref{SDE-X}) is
non-explosive, so that, given any initial condition $x$,
$X_t^0 \in \, ]0,\infty [$, for all $t \geq 0$, with probability 1.

\begin{ass} {\rm
The It\^{o} diffusion $X^0$ defined by (\ref{SDE-X}) is non-explosive.
} \mbox{}\hfill$\Box$ \label{A2} \end{ass}
Feller's test for explosions (see Theorem~5.5.29 in Karatzas and
Shreve~\cite{KS}) provides a necessary and sufficient condition for
this assumption to hold true.
Indeed, if we define
\be
l_{X^0} (x)= \int _c^x \left[ p_{X^0} (x) - p_{X^0} (s) \right] \,
m_{X^0} (ds) ,
\ee
then Assumption~\ref{A2} is satisfied if and only if
$\lim _{x \downarrow 0} l_{X^0} (x) =  \lim_{x \rightarrow \infty}
l_{X^0} (x) = \infty$.

Now, we model the system's controlled dynamics by the SDE
\ben
dX_t = b(X_t) \, dt + dZ_t + \sigma (X_t) \, dW_t , \quad X_0 = x
> 0 , \label{SDE-XZ}
\een
where the controlled process $Z$ is an increasing process.
With each admissible intervention strategy $Z$ (see Definition
\ref{adm} below), we associate the performance criterion
\ben
J_x(Z) = \limsup _{T \rightarrow \infty} \e \left[ \int _0^T
e^{-\Lambda _t^{r (X)}} h(X_t) \, dt - \int _0^T e^{-\Lambda
_t^{r (X)}} k(X_t) \circ dZ_t \right] , \label{J}
\een
where
\ben
\Lambda _t^{r (X)} = \int _0^t r(X_\upsilon) \, d\upsilon ,
\label{Lambda-r}
\een
for some functions $h : \, ]0,\infty [ \, \rightarrow \R$ and
$k, r : \, ]0,\infty [ \, \rightarrow \R _+$.
The integral with respect to $Z$ is defined by
\ben
\int _0^T e^{-\Lambda _t^{r (X)}} k(X_t) \circ dZ_t
= \int _0^T e^{-\Lambda _t^{r (X)}} k(X_t) \, dZ_t^c +
\sum _{0\leq t \leq T} \int _0^{\Delta Z_t} e^{-\Lambda _t^{r (X)}}
k(X_t+s) \, ds , \label{Z-integral}
\een
where $Z^c$ is the continuous part of the increasing process
$Z$.

\begin{de} {\rm
The family ${\mathcal A}$ of all {\em admissible intervention
strategies\/} is the set of all $(\F_t)$-adapted c\`{a}gl\`{a}d
processes $Z$ with increasing sample paths such that $Z_0 = 0$,
(\ref{SDE-XZ}) has a unique non-explosive strong solution, and
\be
\e \left[ \int _0^T e^{-\Lambda _t^{r (X)}} k(X_t) \circ dZ_t
\right] < \infty , \quad \text{for all } T>0 .
\ee
} \label{adm} \end{de}
The objective of our control problem is to maximise $J_x$ over all
admissible strategies.
Accordingly, we define the problem's value function $v$ by
\be
v(x) = \sup _{Z \in {\mathcal A}} J_x(Z) , \quad
\text{for } x>0 .
\ee

%-------------------------------------------------------------
For our optimisation problem to be well-posed and to admit a solution
that conforms with economic intuition, we need to make additional
assumptions.

\begin{ass} {\rm
The discounting factor $r$ is $C^1$, $r'$ is locally Lipschitz, and
there exists a constant $r_0 > 0$ such that
$r(x) \geq r_0$, for all $x>0$.
} \mbox{}\hfill$\Box$ \label{Ar} \end{ass}
Our analysis will also involve the SDE
\ben
dY_t^0 = \mu ( Y_t^0 ) \, dt + \sigma (Y_t^0) \, dW_t , \quad Y_0^0
= x > 0 , \label{SDE-Y}
\een
where
\ben
\mu (x) = b(x) + \sigma (x) \sigma '(x) - \half \sigma ^2
(x) \frac{r'(x)}{r(x)} . \label{mu}
\een
In the presence of Assumptions~\ref{A1} and~\ref{Ar}, this SDE has a
unique strong solution.
Also, a short calculation shows that the scale
function $p_{Y^0}$ of this diffusion satisfies
\ben
p_{Y^0}' (x) := \exp \left( - \int _c^x \frac{2\mu (s)}
{\sigma ^2 (s)} \, ds \right) = \frac{\sigma ^2 (c)}{r(c)}
\frac{r(x)} {\sigma ^2 (x)} p_{X^0} '(x) , \label{p-XY}
\een
where $p_{X^0}$ is the scale function of the diffusion $X^0$ defined
by (\ref{scale}).
We make the following additional assumption.
\begin{ass} {\rm
The It\^{o} diffusion $Y^0$ defined by (\ref{SDE-Y}) is non-explosive.
} \mbox{}\hfill$\Box$ \label{AY} \end{ass}

\begin{rem} {\rm
 It would be of interest to establish conditions under which
 Assumption~\ref{A2} implies Assumption~\ref{AY}.
 We have not been able to address this issue in a comprehensive
 manner.
 However, we note that sufficient conditions can be established
 by appealing to appropriate comparison theorems for solutions to
 SDEs.
 To illustrate this claim, suppose that
\ben
 \sigma (x) \sigma '(x) - \half \sigma ^2 (x) \frac{r'(x)}{r(x)}
 \geq 0, \quad \text{for all } x>0 , \label{rem-r}
\een
and that
\ben
 \mu \text{ is Lipschitz continuous in } [1,\infty[ . \label{rem-mu}
\een
In this case, Assumption~\ref{A2}, (\ref{rem-r}) and the comparison
theorem for solutions to SDE's in Karatzas and Shreve
\cite[Proposition 5.2.18]{KS} imply that $Y^0$ does not explode
at $0$.
Furthermore, (\ref{rem-mu}) implies that $Y^0$ does not explode at
$\infty$.
To see this claim, we argue by contradiction, and we assume that
there exists an $(\F_t)$-stopping time $\tau_\infty$ such that
$\p (\tau_\infty < \infty) > 0$ and $\lim _{t \rightarrow \tau_\infty}
Y_t^0 = \infty$.
Since $Y^0$ does not explode at $0$, there exist $\varepsilon , \zeta
> 0$ such that
\ben
\p \left( \tau_\infty < \infty , \inf _{t \leq \tau_\infty} Y^0_t >
\varepsilon \right) > \zeta . \label{contra}
\een
Given such constants, let $\mu _\varepsilon , \sigma _\varepsilon :
\R \rightarrow \R$ be any Lipschitz continuous functions such that
$\mu _\varepsilon (x) = \mu (x)$ and $\sigma _\varepsilon (x) = \sigma
(x)$, for all $x \geq \varepsilon$, let $Y^\varepsilon$ be the solution
to
\be
dY_t^\varepsilon = \mu _\varepsilon ( Y_t^\varepsilon ) \, dt + \sigma
_\varepsilon (Y_t^\varepsilon) \, dW_t , \quad Y_0^\varepsilon
= x > 0 , 
\ee
and let $\tau_\varepsilon$ be the first hitting time of $]0,\varepsilon]$
by $Y^\varepsilon$.
The Lipschitz continuity of $\mu_\varepsilon$ and $\sigma _\varepsilon$
imply that $Y_t^\varepsilon$ is real-valued, for all $t \geq 0$.
However, this observation provides the required contradiction because
$Y_t^0 = Y_t^\varepsilon$, for all $t \leq \tau_\varepsilon$, so
(\ref{contra}) cannot be true.
\mbox{}\hfill$\Box$ } \end{rem}

Throughout the paper, we denote by $K$ the function defined by
\ben
K(x) = \int _0^x k(s) \, ds , \label{K}
\een
and, given a $C^2$ function $w$, we denote by ${\mathcal L}_X w$
the function given by
\ben
{\mathcal L}_X w(x) = \half \sigma ^2 (x) w'' (x) + b(x) w' (x)
- r(x) w(x) . \label{L}
\een
Also, we define
\ben
Q(x) =  h(x) + {\mathcal L}_X K (x) . \label{Q}
\een

We can now complete the list of assumptions made in the paper.

%-------------------------------------------------------------
\begin{ass} {\rm
The following conditions hold:
\\
(a)
The running payoff function $h$ is $C^1$ and the function $h/r$
is bounded from below. \\
(b)
The running cost function $k$ is $C^2$.
Also, $k(x) \geq 0$, for all $x>0$, and the function $K$ defined by
(\ref{K}) is real-valued. \\
(c)
The problem data is such that
\ben
\rho (x) := \frac{r^2(x) - r(x) b'(x) + r'(x) b(x)}{r(x)} \geq r_0 ,
\quad \text{for all } x>0 , \label{rho}
\een
where the constant $r_0$ is the same as in Assumption~\ref{Ar},
without loss of generality. \\
(d)
There exists a real $x^* \geq 0$ such that
\be
\drq (x) := \frac{r(x) Q'(x) - r'(x) Q(x)}{r(x)} \begin{cases}
\geq 0 , & \text{for } x \leq x^* , \text{ if } x^* > 0 , \\ < 0 ,
& \text{for } x > x^* , \end{cases}
\ee
where the function $Q$ is defined by (\ref{Q}).
Also, if $x^*=0$, then $\lim _{x \downarrow 0} Q(x) / r(x) <
\infty$. \\
(e)
The integrability condition
\ben
\e \left[ \int _0^\infty e^{-\Lambda _t^{r(X^0)}} \left[ \left|
h(X_t^0) \right| + \left| {\mathcal L}_X K(X_t^0) \right| \right]
dt \right] < \infty \label{IC1}
\een
is satisfied, and
\ben
K(x) = - \e \left[ \int _0^\infty e^{-\Lambda _t^{r(X^0)}}
{\mathcal L}_X K(X_t^0) \, dt \right] , \label{IC2}
\een
for every initial condition $x>0$.
} \mbox{}\hfill$\Box$ \label{A3} \end{ass}

%-------------------------------------------------------------
\begin{rem} {\rm
It is worth noting that (\ref{IC2}) is in essence an integrability
condition.
Indeed, an application of It\^{o}'s formula yields
\ben
e^{-\Lambda _T^{r(X^0)}} K(X_{T+}^0) = K(x) + \int _0^T e^{-\Lambda
_t^{r(X^0)}} {\mathcal L}_X K(X_t^0) \, dt +  \int _0^T e^{-\Lambda
_t^{r(X^0)}} \sigma (X_t^0) k(X_t^0) \, dW_t . \label{IC2Ito}
\een
If we assume that the stochastic integral appearing in this identity
is a martingale and that the so-called transversality condition
\be
\liminf _{T \rightarrow \infty} \e \left[ e^{-\Lambda _T^{r(X^0)}}
K(X_{T+}^0) \right] = 0
\ee
holds, then we can take expectations in (\ref{IC2Ito}) and then pass
to the limit $T \rightarrow \infty$ to obtain (\ref{IC2}).
} \mbox{}\hfill$\Box$ \end{rem}

%-------------------------------------------------------------
\begin{rem} {\rm
For future reference, we observe that the integrability condition
(\ref{IC1}) implies that
\be
\e \left[ \int _0^\infty e^{-\Lambda _t^{r(X^0)}} \left| Q(X_t^0)
\right| dt \right] < \infty . %\label{IC3}
\ee
Furthermore, if we define
\ben
R_{X^0,h}(x) = \e \left[ \int _0^\infty e^{-\Lambda _t^{r (X^0)}}
h(X_t^0) \, dt \right] , \label{Rh-def}
\een
then
\ben
R_{X^0,h}(x) - K(x) = \e \left[ \int _0^\infty e^{-\Lambda _t
^{r (X^0)}} Q(X_t^0) \, dt \right] =: R_{X^0,Q}(x) , \label{Rh-K/Q}
\een
thanks to (\ref{IC2}).
} \mbox{}\hfill$\Box$ \label{rem:RhQ}\end{rem}

%-------------------------------------------------------------
The following result is concerned with the well-posedness of our
optimisation problem.

\begin{lem}
Suppose that Assumptions~\ref{A1}, \ref{A2}, \ref{Ar}, \ref{A3}.(a),
\ref{A3}.(b) and \ref{A3}.(d) hold true, and fix any initial condition
$x>0$.
The performance index $J_x(Z)$ is well-defined for every admissible
intervention strategy $Z \in {\mathcal A}$, and $v(x)<\infty$.
\end{lem}
{\bf Proof.}
%===========
Fix any initial condition $x>0$ and any admissible intervention process
$Z \in {\mathcal A}$.
Using It\^{o}'s formula, the fact that $\Delta X_t = \Delta Z_t$ and
the definition (\ref{K}) of the function $K$, we obtain
\begin{align}
e^{-\Lambda _T^{r (X)}} K(X_{T+}) = \mbox{} & K(x) + \int _0^T
e^{-\Lambda _t^{r (X)}} {\mathcal L}_X K(X_t) \, dt + \int _0^T
e^{-\Lambda _t^{r (X)}} k(X_t) \, dZ_t  \nonumber \\
& + \sum _{0 \leq t \leq T} e^{-\Lambda _t^{r (X)}} \left[ K(X_{t+}) -
K(X_t) - k(X_t) \Delta X_t \right] + M_T \nonumber \\
= \mbox{} & K(x) + \int _0^T e^{-\Lambda _t^{r (X)}} {\mathcal L}_X
K(X_t) \, dt + \int _0^T e^{-\Lambda _t^{r (X)}} k(X_t) \, dZ_t^c
\nonumber \\
& + \sum_{0 \leq t \leq T} e^{-\Lambda _t^{r (X)}} \int _{X_t}^{X_{t+}}
k(s) \, ds + M_T , \nonumber
\end{align}
where $Z^c$ is the continuous part of the increasing process $Z$ and
\be
M_T = \int_0^T e^{-\Lambda _t^{r (X)}} \sigma(X_t) k(X_t) \, dW_t.
\ee
In view of (\ref{Z-integral}) and the definition (\ref{Q}) of $Q$,
this calculation implies that
\begin{align}
\int _0^T e^{-\Lambda _t^{r (X)}} h(X_t) \, dt - \int _0^T
e^{-\Lambda _t^{r (X)}} k(X_t) \circ dZ_t = \mbox{} & K(x) -
e^{-\Lambda _T^{r (X)}} K(X_{T+}) \nonumber \\
& + \int _0^T e^{-\Lambda _t^{r (X)}} Q(X_t) \, dt + M_T .
\label{Well-Ito}
\end{align}

Now, combining the assumption that $r(x) \geq r_0 >0$, for all $x>0$,
with Assumption~\ref{A3}.(d), we can see that the function $Q/r$ is
bounded from above.
It follows that
\ben
\int _0^T e^{-\Lambda _t^{r (X)}} Q(X_t) \, dt \leq
\sup _{s>0} \frac{Q(s)}{r(s)} \int_0^\infty e^{-\Lambda _t^{r (X)}}
r(X_t) \, dt = \sup _{s>0} \frac{Q(s)}{r(s)} < \infty
. \label{Q-int-bound}
\een
Also, the assumption that $h/r$ is bounded from below implies that
\be
\int _0^T e^{-\Lambda _t^{r (X)}} h(X_t) \, dt \geq \inf _{s>0}
\frac{h(s)}{r(s)} \int_0^\infty e^{-\Lambda _t^{r (X)}} r(X_t) \, dt
\geq \inf _{s>0} \frac{h(s)}{r(s)} > - \infty .
\ee
In light of these inequalities, the fact that $K$ is positive
and real-valued, and (\ref{Well-Ito}), we can see that
\be
\inf _{t \leq T} M_t \geq \inf _{s>0} \frac{h(s)}{r(s)} - \sup
_{s>0} \frac{Q(s)}{r(s)} - K(x) - \int _0^T e^{-\Lambda
_t^{r (X)}} k(X_t) \circ dZ_t .
\ee
The admissibility of $Z$ implies that the right-hand side of this
inequality belongs to ${\mathcal L}^1$, for all $T>0$.
It follows that $M$ is a supermartingale.
Combining this observation with (\ref{Q-int-bound}), we can see that
(\ref{Well-Ito}) implies that
\be
\e \left[ \int _0^T e^{-\Lambda _t^{r (X)}} h(X_t) \, dt - \int
_0^T e^{-\Lambda _t^{r (X)}} k(X_t) \circ dZ_t \right] \leq
\sup _{s>0} \frac{Q(s)} {r(s)} + K(x) < \infty .
\ee
However, these inequalities establish the claims made.
\mbox{}\hfill$\Box$

%=============================================================
\section{The solution to a second order linear ODE}
\label{sec:ODE}

In this section, we review some properties of the solution to
the ODE
\ben
{\mathcal L}_X w(x) +G(x) \equiv \half \sigma ^2 (x) w''(x) + b(x)
w'(x) -r(x) w(x) + G(x) = 0 , \label{ODE-X}
\een
that is associated with the It\^{o} diffusion (\ref{SDE-X}), and
we establish some new results that are concerned with appropriate
analytic and probabilistic expressions for the derivative $w'$ of a
solution $w$ to (\ref{ODE-X}).
Indeed, if $G$ is absolutely continuous, then differentiating the
ODE (\ref{ODE-X}) yields
\be
\half \sigma ^2 (x) w'''(x) + \left[ b(x) + \sigma (x)
\sigma '(x) \right] w''(x) - \left[ r(x) -b'(x) \right]
w'(x) - r'(x) w(x) + G'(x) = 0 .
\ee
Using (\ref{ODE-X}) once again to eliminate $w(x)$ from
this equation, we can see that $w'$ solves the ODE
\ben
{\mathcal L}_Y u (x) + \drg (x) := \half \sigma ^2 (x) u''
(x) + \mu (x) u' (x) - \rho (x) u(x) + \drg (x) = 0 ,
\label{ODE-Y}
\een
where $\mu$ is defined by (\ref{mu}), $\rho$ is defined by
(\ref{rho}), and
\ben
\drg (x) := \frac{r(x) G'(x) - r'(x) G(x)}{r(x)} = r(x) \frac{d}
{dx} \left( \frac{G(x)}{r(x)} \right) .
\label{drg}
\een
It follows that $w'$ satisfies a second order linear ODE that is
similar to (\ref{ODE-X}) and is associated with the SDE (\ref{SDE-Y}).

In the presence of Assumptions~\ref{A1}, \ref{A2} and
\ref{Ar}, the general solution to the homogeneous ODE
${\mathcal L}_X w(x) =0$, which is associated with (\ref{ODE-X}), 
is given by
\be
w(x) = A \phi (x) + B \psi (x) ,
\ee
for some constants $A, B \in \R$, where $\phi$ and $\psi$ are
$C^2$ functions such that
\begin{gather}
0 < \phi (x) \quad \text{and} \quad \phi ' (x) < 0 , \quad
\text{for all } x>0 , \label{phi-psi-prop1} \\
0 < \psi (x) \quad \text{and} \quad \psi ' (x) > 0 , \quad
\text{for all } x>0 , \label{phi-psi-prop2} \\
\intertext{and}
\lim _{x \downarrow 0} \phi (x) = \lim _{x \rightarrow \infty}
\psi (x) = \infty . \label{phi-psi-prop3}
\end{gather}
These functions are unique, modulo multiplicative constants.
To simplify the notation we assume, without loss of generality,
that
\ben
\phi (c) = \psi (c) = 1 , \label{phi-psi-c}
\een
where $c>0$ is the same constant as the one that we used in the
definition (\ref{scale}) of the scale function $p_{X^0}$.
Also, these functions satisfy
\ben
\phi (x) \psi '(x) - \phi ' (x) \psi (x) = C p_{X^0}'(x) ,
\label{phi-psi-scale}
\een
where
\ben
C := \left[ \psi '(c) - \phi '(c) \right] > 0 . \label{C}
\een
Furthermore, we can use the fact that $\phi$ and $\psi$ satisfy the
ODE ${\mathcal L}_X w(x) = 0$ to verify that
\begin{align}
\phi '' (x) \psi '(x) - \phi '(x) \psi '' (x) & = \frac{2r(x)}
{\sigma ^2 (x)} \left[ \phi (x) \psi '(x) - \phi '(x) \psi (x)
\right] \nonumber \\
& = \frac{2 C r(x)}{\sigma ^2(x)} p_{X^0}'(x) . \label{phi-psi'-Wr}
\end{align}

Similarly, Assumptions~\ref{A1}, \ref{Ar}, \ref{AY} and \ref{A3}.(c)
guarantee that the general solution to the homogeneous ODE ${\mathcal
L}_Y u(x) = 0$, which is associated with (\ref{ODE-Y}), is given by
\be
u(x) = \tilde{A} \tilde{\phi}(x) + \tilde{B} \tilde{\psi}
(x) , \ee
for some constants $\tilde{A} , \tilde{B} \in \R$, where
$\tilde{\phi}$ and $\tilde{\psi}$ are $C^2$ functions
satisfying
\begin{gather}
0 < \tilde{\phi} (x) \quad \text{and} \quad \tilde{\phi}'
(x) < 0 , \quad \text{for all } x>0 ,
\label{phi-psi-t-prop1} \\
0 < \tilde{\psi} (x) \quad \text{and} \quad \tilde{\psi}'
(x) > 0 , \quad \text{for all } x>0 , \label{phi-psi-t-prop2} \\
\intertext{and}
\lim _{x \downarrow 0} \tilde{\phi} (x) = \lim
_{x \rightarrow \infty} \tilde{\psi} (x) = \infty ,
\label{phi-psi-t-prop3}
\end{gather}
Once again, we assume, without loss of generality, that
\ben
\tilde{\phi} (c) = \tilde{\psi} (c) = 1 , \label{phi-psi-c-til}
\een
for notational simplicity.
Also, we note that
\be
\tilde{\phi} (x) \tilde{\psi}' (x) - \tilde{\phi}' (x) \tilde{\psi}
(x) = \tilde{C} p_{Y^0}' (x) ,
\ee
where
\ben
\tilde{C} := \tilde{\psi}' (c) - \tilde{\phi}' (c) > 0 . \label{C-til}
\een

The following example shows that Assumption~\ref{A3}.(c) is
indispensable if we want the functions $\tilde{\phi}$ and
$\tilde{\psi}$ to have the properties
(\ref{phi-psi-t-prop1})--(\ref{phi-psi-t-prop3}), which are
essential for constructing the solutions to many one-dimensional
stochastic control problems, including the one that we study in
this paper.

%-------------------------------------------------------------
\begin{ex} {\rm
Suppose that $b(x)=2x$, $\sigma(x)=\sqrt{2} x$ and $r(x) =
\frac{3}{4}$.
In this case, $\rho (x) = - \frac{5}{4}$ and
Assumption~\ref{A3}.(c) is not satisfied.
A simple calculation reveals that the functions
\be
\phi (x) = x^{-\frac{3}{2}} , \quad \psi (x) = x^{\frac{1}{2}}
\ee
span the solution space of the ODE ${\mathcal L}_X w(x) = 0$,
and that the functions
\be
\tilde{\phi} (x) = x^{-\frac{5}{2}} , \quad \tilde{\psi} (x) =
x^{-\frac{1}{2}}
\ee
span the solution space of the ODE ${\mathcal L}_Y u(x) = 0$.
Clearly, $\tilde{\psi}$ does not satisfy (\ref{phi-psi-t-prop2}) and
(\ref{phi-psi-t-prop3}).
} \end{ex}

The next result expresses the functions $\tilde{\phi}$ and
$\tilde{\psi}$ in terms of the functions $\phi$ and $\psi$.

%-------------------------------------------------------------
\begin{prop}
Suppose that Assumptions~\ref{A1}, \ref{A2}, \ref{Ar}, \ref{AY}
and~\ref{A3}.(c) hold.
If $\phi$, $\psi$ and $\tilde{\phi}$, $\tilde{\psi}$ are the
functions satisfying (\ref{phi-psi-prop1})--(\ref{phi-psi-c})
and (\ref{phi-psi-t-prop1})--(\ref{phi-psi-c-til}), and
spanning the solution space of the homogeneous ODEs
${\mathcal L}_X w(x) = 0$ and ${\mathcal L}_Y u(x) = 0$
associated with (\ref{ODE-X}) and (\ref{ODE-Y}), respectively,
then
\ben
\tilde{\phi} (x) = \frac{1}{\phi' (c)} \phi' (x) \quad \text{and} \quad
\tilde{\psi} (x) = \frac{1}{\psi' (c)} \psi' (x) . \label{phi-psi/til}
\een
Also, both of $\phi$ and $\psi$ are convex, and, if $C$ and
$\tilde{C}$ are the constants defined by (\ref{C}) and (\ref{C-til}),
respectively, then
\ben
\tilde{C} = - \frac{2 r(c)}{\sigma ^2 (c)} \frac{1}{\phi' (c) \psi'
(c)} C . \label{C-Ctil}
\een
\label{prop:phi-psi's} \end{prop}
{\bf Proof.}
%===========
We first show that $\tilde{\psi} = C_\psi \psi '$, for some constant
$C_\psi > 0$.
To this end, we define the function $\hat{\psi}$ by
\ben
\hat{\psi}(x) = \int _1^x \tilde{\psi} (s) \, ds . \label{hat-psi}
\een
Given an absolutely continuous function $f : \, ]0,\infty [
\, \rightarrow \R$ with compact support, we can use the integration
by parts formula to calculate
\be
\int_0^\infty f'(s) {\mathcal L}_X \hat{\psi} (s) \, ds =
- \int_0^\infty f(s) \left[ {\mathcal L}_Y \hat{\psi}' (s)
+ \frac{r'(s)}{r(s)}{\mathcal L}_X \hat{\psi} (s) \right] ds
\ee
Since $\hat{\psi}' = \tilde{\psi}$ satisfies the ODE ${\mathcal L}
_Y u(x) =0$, it follows that
\ben
\int_0^\infty \frac{r(s) f'(s) + r'(s) f(s)}{r(s)} {\mathcal L}
_X \hat{\psi}(s) \, ds =0 . \label{int-by-parts}
\een
Now, fix any $a>0$ and any $z \in \, ]0,a[$, and define
\be
q_{a,z} (x) = \begin{cases} 1 , & \text{if } x \in
[a-z, a] , \\ -1 , & \text{if } x \in \,]a, a+z] ,
\\ 0 , & \text{otherwise}. \end{cases}
\ee
If we choose
\be
f(x) = \frac{1}{r(x)} \int _0^x q_{a,z} (s) \, ds ,
\ee
then $f$ is absolutely continuous with compact support,
and (\ref{int-by-parts}) yields
\be
\int _{a-z}^a \frac{{\mathcal L} _X \hat{\psi}(s)}
{r(s)} \, ds = \int _a^{a+z} \frac{{\mathcal L} _X
\hat{\psi}(s)} {r(s)} \, ds .
\ee
Since $z \in \, ]0,a[$ has been arbitrary, we can differentiate
this expression with respect to $z$ to obtain
\be
\frac{{\mathcal L} _X \hat{\psi}(a-z)} {r(a-z)}  =
\frac{{\mathcal L} _X \hat{\psi}(a+z)} {r(a+z)}  , \quad
\text{for all } z \in \, ]0,a[ .
\ee
It follows that the function ${\mathcal L} _X \hat{\psi} / r$
has even symmetry around the point $a$.
However, since $a>0$ has been arbitrary, this can be true
only if ${\mathcal L} _X \hat{\psi} / r$ is constant.
Therefore, there exists a constant $C_1 \in \R$ such that
${\mathcal L} _X \left( \hat{\psi}(x) + C_1 \right)  =0$.
If we combine this observation with the facts that
\be
\hat{\psi} \text{ is strictly increasing,} \quad \lim
_{x \downarrow 0} \hat{\psi}(x) \in \R \quad \text{and} \quad
\lim _{x \rightarrow \infty} \hat{\psi}(x) = \infty ,
\ee
which follow from the definition (\ref{hat-psi}) of
$\hat{\psi}$ and the properties
(\ref{phi-psi-t-prop2})--(\ref{phi-psi-t-prop3}) of the function
$\tilde{\psi}$, and the properties
(\ref{phi-psi-prop1})--(\ref{phi-psi-prop3}) of the functions
$\phi$, $\psi$ that span the solution space of ${\mathcal L}
_X w(x)$, we can conclude that $\hat{\psi} + C_1 = C_\psi \psi$,
for some constant $C_\psi > 0$.
However, this conclusion establishes the identity $\tilde{\psi}
= C_\psi \psi'$.

Now, since $\phi'$ and $\psi' = C_\psi^{-1} \tilde{\psi}$ are
independent solutions to the ODE ${\mathcal L}_Y u(x) = 0$ and
${\mathcal L}_Y \tilde{\phi}(x) = 0$, there exist constants
$C_\phi$ and $\Gamma$ such that $\tilde{\phi} = - C_\phi
\phi ' + \Gamma \tilde{\psi}$.
However, the limits
\be
\lim _{x \rightarrow \infty} \tilde{\phi}(x) \in [0,\infty[ ,
\quad \lim _{x \rightarrow \infty} \phi ' (x) = 0  \quad \text{and}
\quad \lim _{x \rightarrow \infty} \tilde{\psi}(x) = \infty
\ee
imply that $\Gamma = 0$, which proves that $\tilde{\phi} = -
C_\phi \phi '$, for some constant $C_\phi >0$.

Finally, we note that
\be
- C_\phi = \frac{1}{\phi' (c)} \quad \text{and} \quad C_\psi =
\frac{1}{\psi' (c)}
\ee
are the only choices for the constants $C_\phi$ and $C_\psi$ that are
compatible with (\ref{phi-psi-c-til}).
Also, (\ref{C-Ctil}) follows by a simple calculation involving the
definitions of the constants $C$ and $\tilde{C}$, (\ref{phi-psi/til})
and (\ref{phi-psi'-Wr}).
\mbox{}\hfill$\Box$
\vspace{5mm}

%-------------------------------------------------------------
To proceed further, we consider the solution $X^0$ to the SDE
(\ref{SDE-X}) and we define
\be
\Lambda _t^{r (X^0)} = \int _0^t r(X_\upsilon^0) \, d\upsilon .
\ee
We recall that, if $G$ is a measurable function,
then
\ben
\e \left[ \int _0^\infty e^{-\Lambda _t^{r (X^0)}} \left| G(X_t^0)
\right| dt \right] < \infty \label{IC-X1}
\een
if and only if
\ben
\int_0^x \frac{|G(s)| \psi (s)}{\sigma^2(s) p_{X^0}' (s)} \, ds
+ \int _x^\infty \frac{|G(s)| \phi (s)}{\sigma^2(s) p_{X^0}' (s)}
\, ds < \infty . \label{IC-X2}
\een
In the presence of these equivalent integrability conditions,
the function $R_{X^0,G}$ defined by
\ben
R_{X^0,G}(x) = \e \left[ \int _0^\infty e^{-\Lambda _t^{r (X^0)}}
G(X_t^0) \, dt \right] , \label{RG-def}
\een
admits the analytic expression
\ben
R_{X^0,G} (x) = \frac{2}{C} \phi (x) \int_0^x \frac{G(s) \psi
(s)}{\sigma^2(s) p_{X^0}' (s)} \, ds + \frac{2}{C} \psi (x) \int
_x^\infty \frac{G(s) \phi (s)} {\sigma^2(s) p_{X^0}' (s)} \, ds ,
\label{RG}
\een
where $C>0$ is the constant defined by (\ref{C}), and is a
special solution to (\ref{ODE-X}).

At this point, we establish the following technical result that we
will need.

%-------------------------------------------------------------
\begin{lem}
Suppose that Assumptions~\ref{A1}, \ref{A2} and \ref{Ar} hold.
If $G$ is a function satisfying (\ref{IC-X1}) and (\ref{IC-X2}), 
then
\ben
\liminf _{x \downarrow 0} \frac{|G(x)| \psi' (x)}{r(x)
p_{X^0}'(x)} = \liminf _{x \rightarrow \infty} \frac{|G(x)|
\phi' (x)}{r(x) p_{X^0}'(x)}  =0 . \label{int-by-parts-lims}
\een
\label{liminf-lemma} \end{lem}
{\bf Proof.}
%===========
In view of (\ref{phi-psi-scale}), we can see that
\be
0 < \frac{\phi (x) \psi ' (x)}{C p_{X^0}'(x)} < 1 \quad
\text{and} \quad 0 < - \frac{\phi ' (x) \psi (x)}{C p_{X^0}'(x)}
< 1 ,
\ee
which, combined with (\ref{phi-psi-prop3}) implies that
\ben
\lim _{x \downarrow 0} \frac{\psi'(x)}{p_{X^0}'(x)} = \lim _{x
\rightarrow \infty} \frac{\phi ' (x)}{p_{X^0}' (x)} = 0
. \label{IC-X3}
\een
Also, the calculation
\be
\frac{d}{dx} \left( \frac{1}{p_{X^0}'(x)} \right) = \frac{2b(x)}
{\sigma ^2 (x) p_{X^0}'(x)}
\ee
and the fact that $\phi$ satisfies the ODE ${\mathcal L}_X w (x) =
0$, imply that
\begin{align}
\frac{d}{dx} \left( \frac{\phi ' (x)}{p_{X^0}'(x)} \right) & =
\frac{2}{\sigma ^2 (x) p_{X^0}'(x)} \left[ \half \sigma ^2 (x)
\phi '' (x) + b(x) \phi ' (x) \right] \nonumber \\
& = \frac{2 r(x) \phi (x)}{\sigma ^2 (x) p_{X^0}'(x)} .
\label{IC-X4}
\end{align}
Similarly, we can show that
\ben
\frac{d}{dx} \left( \frac{\psi ' (x)}{p_{X^0}'(x)} \right) =
\frac{2 r(x) \psi (x)}{\sigma ^2 (x) p_{X^0}'(x)} . \label{IC-X5}
\een

Now, we consider any sequence $(x_n)$ such that
\be
0 < x_n < \frac{1}{n} \quad \text{and} \quad \frac{|G(x_n)|}{r(x_n)}
\leq \inf _{x \in \, \left] 0 , \frac{1}{n} \right[}
\frac{|G(x)|}{r(x)} + 1 .
\ee
Using (\ref{IC-X3}) and (\ref{IC-X4}), we calculate
\begin{align}
\int _0^{x_n} \frac{2|G(s)| \psi(s)}{\sigma ^2(s)p_{X^0} '(s)} \, ds
& \geq \left( \frac{|G(x_n)|}{r(x_n)} - 1 \right) \int
_0^{x_n} d\left( \frac{\psi'(s)}{p_{X^0} '(s)} \right) \nonumber \\
& = \left( \frac{|G(x_n)|}{r(x_n)} - 1 \right)
\frac{\psi'(x_n)}{p_{X^0} '(x_n)} . \nonumber
\end{align}
In view of (\ref{IC-X3}) and the fact that $\lim _{n \rightarrow
\infty} \int _0^{x_n} \frac{2|G(s)| \psi(s)}{\sigma ^2(s)p_{X^0} '(s)}
\, ds = 0$, we can pass to the limit $n\rightarrow\infty$ in this
inequality to obtain
\be
\lim _{n \rightarrow \infty} \frac{|G(x_n)|\psi '(x_n)}{r(x_n)
p_{X^0} '(x_n)}=0 ,
\ee
which proves that
$\liminf _{x \downarrow 0} \frac{|G(x)| \psi' (x)}{r(x)
p_{X^0}'(x)}=0$.
Using similar arguments, we can also show that
$\liminf _{x \rightarrow \infty} \frac{|G(x)|\phi' (x)}{r(x)
p_{X^0}'(x)} =0$.
\mbox{}\hfill$\Box$
\vspace{5mm}

%-------------------------------------------------------------
\begin{rem} {\rm
It is worth noting that the conclusions of the preceding result cannot
be strengthened.
To see this, suppose that $b(x)=-x$, $\sigma(x)=\sqrt{2} x$ and $r(x)=3$,
so that
\be
\phi(x)=1/x , \quad \psi(x)=x^3 \quad \text{and} \quad p_{X^0} '(x)=x .
\ee
Also, let $G$ be the positive function given by
\be
G(x) = \sum _{n=1}^{\infty} G^{(n)} (x) ,
\ee
where $G^{(n)}$ is the tent-like function defined, for $n\geq 1$ by
\be
G^{(n)} (x) =
\begin{cases}
n^4 x + n -n^3, & \text{for } x \in \left[\frac{1}{n} - \frac{1}{n^3},
\frac{1}{n} \right[ , \\
-n^4 x + n +n^3, & \text{for } x \in \left[\frac{1}{n}, \frac{1}{n} +
\frac{1}{n^3} \right[ , \\
0, & \text{ otherwise} .
\end{cases}
\ee
Given any $x>0$, we calculate
\begin{align}
\int _0^x \frac{|G(s)| \psi(s)}{\sigma^2 (s) p_{X^0} '(s)} ds & \leq
\int _0^2 G(s) ds \nonumber \\
& = \sum _{n=1} ^{\infty} \frac{1}{2} n \left[ \frac{1}{n} +
\frac{1}{n^3} - \left( \frac{1}{n} - \frac{1}{n^3} \right) \right]
\nonumber \\
& < \infty , \nonumber
\end{align}
and we can immediately see that
\be
\int _x^\infty \frac{|G(s)| \phi(s)}{\sigma ^2 (s) p_{X^0} '(s)} ds <
\infty .
\ee
It follows that the assumptions of Lemma~\ref{liminf-lemma} are
satisfied.
However,
\be
\limsup _{x \downarrow 0} \frac{|G(x)| \psi '(x)}{r(x) p_{X^0} '(x)}
\geq \lim _{n \rightarrow \infty} \frac{1}{n} G^{(n)}
\left(\frac{1}{n}\right) > 0 .
\ee
\mbox{}\hfill$\Box$ } \label{remark-triangle-height-n} \end{rem}
\vspace{5mm}

%-------------------------------------------------------------
In what follows, we assume that $G$ is absolutely continuous.
Also, we consider the solution $Y^0$ to the SDE (\ref{SDE-Y}) and we
define
\be
\Lambda _t^{\rho (Y^0)}= \int _0^t \rho (Y_\upsilon^0 ) \, d\upsilon ,
\ee
If $\drg$ is defined by (\ref{drg}), then the standard theory
discussed above Lemma~\ref{liminf-lemma} implies that
\ben
\e \left[ \int _0^\infty e^{-\Lambda _t^{\rho (Y^0)}} \left| \drg
(Y_t^0) \right| dt \right] < \infty \label{IC-Y1}
\een
if and only if
\ben
\int_0^x \frac{|\drg (s)| \tilde{\psi} (s)}{\sigma^2(s) p_{Y^0}' (s)}
\, ds + \int _x^\infty \frac{|\drg (s)| \tilde{\phi} (s)}{\sigma
^2(s) p_{Y^0}' (s)} \, ds < \infty . \label{IC-Y2}
\een
In light of (\ref{phi-psi/til}) in Proposition~\ref{prop:phi-psi's}
and (\ref{p-XY}), we can see that (\ref{IC-Y2}) is true if
and only if
\ben
\int_0^x \frac{|\drg (s)| \psi' (s)}{r(s) p_{X^0}' (s)} \, ds
+ \int _x^\infty \frac{|\drg (s)| [-\phi '] (s)}{r(s) p_{X^0}' (s)}
\, ds < \infty . \label{IC-Y3}
\een
Furthermore, if these equivalent conditions hold true, then the
function $R_{Y^0, \drg}$ defined by
\be
R_{Y^0, \drg} (x) = \e \left[ \int _0^\infty e^{-\Lambda _t
^{\rho (Y^0)}} \drg (Y_t^0) \, dt \right]
\ee
admits the analytic expressions
\begin{align}
R_{Y^0,\drg} (x) & = \frac{2}{\tilde{C}} \tilde{\phi} (x) \int_0^x
\frac{(\drg) (s) \tilde{\psi} (s)}{\sigma^2 (s) p_{Y^0}'(s)} \, ds +
\frac{2}{\tilde{C}} \tilde{\psi} (x) \int _x^\infty \frac{(\drg) (s)
\tilde{\phi} (s)} {\sigma^2 (s) p_{Y^0}'(s)} \, ds \nonumber \\
& = \frac{1}{C} [-\phi'] (x) \int_0^x \frac{(\drg) (s) \psi' (s)}
{r(s) p_{X^0}'(s)} \, ds + \frac{1}{C} \psi' (x) \int _x^\infty
\frac{(\drg) (s) [-\phi'] (s)} {r(s) p_{X^0}'(s)} \, ds ,
\label{Rdrg}
\end{align}
where $C$ and $\tilde{C}$ are defined by (\ref{C}) and (\ref{C-til}),
respectively, and provides a solution to the ODE ${\mathcal L}_Y u (x)
+ \drg (x) =0$.
Once again, the second equality here follows from
Proposition~\ref{prop:phi-psi's} and (\ref{p-XY}).

%-------------------------------------------------------------
\begin{prop}
Suppose that Assumptions~\ref{A1}, \ref{A2}, \ref{Ar}, \ref{AY} and
\ref{A3}.(c) hold true, and let $G$ be a function satisfying
(\ref{IC-X1}) and (\ref{IC-X2}).
Also, suppose that there exists a constant $\varepsilon>0$ such that
\begin{gather}
\text{either } \drg (s) \leq 0 , \text{ for all } s \leq \varepsilon ,
\text{ or } \drg (s) \geq 0 , \text{ for all } s \leq \varepsilon
\nonumber \\
\intertext{and}
\text{either } \drg (s)\leq 0 , \text{ for all } s \geq 1/\varepsilon ,
\text{ or } \drg (s) \geq 0 , \text{ for all } s \geq 1/\varepsilon .
\nonumber
\end{gather}
Under these conditions, the integrability condition (\ref{IC-Y3}) holds
true, the function $R_{Y^0,\drg}$ is well-defined and real-valued, and
\be
R_{X^0,G}' (x) = R_{Y^0,\drg} (x) , \quad \text{for all } x>0 .
\ee
\label{R-(X,G)'=R-(Y,H)} \end{prop}
{\bf Proof.}
%===========
In view of Lemma~\ref{liminf-lemma}, (\ref{IC-X4}), (\ref{IC-X5})
and the relationship (\ref{drg}) of the functions $G$ and $\drg$,
we can use the integration by parts formula to obtain
\begin{gather}
\int _0^x \frac{(\drg) (s) \psi' (s)}{r(s) p_{X^0}'(s)} \, ds =
\frac{G(x)}{r(x)} \frac{\psi' (x)}{p_{X^0}'(x)} - 2 \int _0^x
\frac{G(s) \psi (s)}{\sigma^2 (s) p_{X^0}'(s)} \, ds \nonumber \\
\intertext{and}
\int _x^\infty \frac{(\drg) (s) \left[ -\phi ' \right] (s)} {r(s)
p_{X^0}'(s)} \, ds = \frac{G(x)}{r(x)} \frac{\phi' (x)}{p_{X^0}'(x)} +
2 \int _x^\infty \frac{G(s) \phi (s)}{\sigma^2 (s) p_{X^0}'(s)} \, ds
. \nonumber
\end{gather}
However, combining these expressions and the expressions (\ref{RG})
and (\ref{Rdrg}) for $R_{X^0,G}$ and $R_{Y^0,\drg}$ with the
assumptions on $\drg$ and its local integrability we can see that all
of the statements made hold true.
\mbox{}\hfill$\Box$
\vspace{5mm}

%-------------------------------------------------------------
\begin{rem} {\rm
If we remove the assumption that $\drg (x)$ has constant sign for
all $x$ sufficiently small and for all $x$ sufficiently large, then
the conclusions of the result above, (\ref{IC-Y2}) in particular,
do not necessarily hold.
To see this, suppose that $b$, $\sigma$ and $r$ are as in
Remark~\ref{remark-triangle-height-n}, and let $G$ be the function
given by
\be
G(x) = \sum _{n=1}^{\infty} G^{(n)} (x) ,
\ee
where $G^{(n)}$ is the tent-like function defined by
\be
G^{(n)} (x) =
\begin{cases}
2n (2n+1)x - 2n, & \text{for } x \in \left[ \frac{1}{2n+1} ,
\frac{1}{2n} \right] , \\
-2n (2n-1)x + 2n, & \text{for } x \in \left[ \frac{1}{2n} ,
\frac{1}{2n-1} \right] , \\
0, & \text{ otherwise} .
\end{cases}
\ee
Plainly, $G$ satisfies (\ref{IC-X1}) and (\ref{IC-X2}) because
it is bounded.
However, the calculation
\begin{align}
\int _0^x \frac{|\drg (s)| \psi '(s)}{r(s) p_{X^0} '(s)}ds & =
\int_0^x s |\drg (s)| \, ds \nonumber \\
& = \sum _{n=1}^\infty 2n (2n+1) \int _{\frac{1}{2n+1}}^{\frac{1}{2n}}
s \, ds + \sum _{n=1}^\infty 2n (2n-1) \int
_{\frac{1}{2n}}^{\frac{1}{2n-1}} s \, ds \nonumber \\
& = \sum _{n=1}^\infty \frac{4n+1}{4n(2n+1)} + \sum _{n=1}^\infty
\frac{4n-1}{4n(2n-1)} \nonumber \\
& = \infty \nonumber
\end{align}
shows that (\ref{IC-Y3}) is not satisfied.
\mbox{}\hfill$\Box$ }
\end{rem}
\vspace{5mm}

%-------------------------------------------------------------
The optimal strategy of the control problem that we solve in the next
section reflects the state process in a given point $a>0$ in the
positive direction.
Such a strategy involves the construction of a continuous increasing
process $Z^a$ such that, if $X^a$ is the associated solution to the
SDE (\ref{SDE-XZ}) with $x \geq a$, then
\ben
X_t^a \geq a \quad \text{and} \quad Z_t^a = \int _0^t {\bf 1} _{\{
X_s^a = a \}} \, dZ_s^a , \quad \text{for all } t \geq 0 . \label{XZ-a}
\een
Such a construction is standard and can be found in El Karoui and
Chaleyat-Maurel~\cite{EC} (see also Schmidt~\cite{Sc}).

%-------------------------------------------------------------
\begin{lem}
Suppose that Assumptions~\ref{A1}, \ref{A2}, \ref{Ar}, \ref{AY} and
\ref{A3}.(c) hold.
Given a real number $a>0$, consider the continuous increasing process
$Z^a$ and the solution $X^a$ to the SDE (\ref{SDE-XZ}) with initial
condition $x \geq a$ that satisfy (\ref{XZ-a}).
If $G : \, ]0,\infty [ \, \rightarrow \R$ is a measurable function
satisfying (\ref{IC-X1}) and (\ref{IC-X2}), then
\begin{gather}
\e \left[ \int _0^\infty e^{-\Lambda _t^{r (X^a)}} \left| G (X_t^a)
\right| dt \right] < \infty , \label{IC-Xa} \\
\intertext{and}
U_{X^a,G} (x) := - \frac{R_{X^0,G}' (a)}{\phi ' (a)} \phi (x) +
R_{X^0,G} (x) = \e \left[ \int _0^\infty e^{-\Lambda _t^{r (X^a)}} G
(X_t^a) \, dt \right] , \label{UG}
\end{gather}
where $R_{X^0,G}$ is defined by (\ref{RG}).
Furthermore,
\ben
\lim _{n \rightarrow \infty} \e \left[ e^{-\Lambda _{T_n^a}^{r (X^a)}}
\right] R_{X^0,G} (n) = 0 , \label{TVC-refl}
\een
where $T_n^a$ is the first hitting time of $\{ n \}$ defined by
\ben
T_n^a = \inf \{ t \geq 0 \mid \ X_t^a = n \} . \label{Tan}
\een
\label{refl-lemma} \end{lem}
{\bf Proof.}
%===========
In view of the continuity of $Z^a$ and the fact that it increases on
the set $\{ X_t^a = a \}$ (see (\ref{XZ-a})), we can see that It\^{o}'s
formula and the identities ${\mathcal L}_X U_{X^a,G} (x) + G (x) = 0$
and $U_{X^a,G}' (a) = 0$ imply that
\begin{align}
e^{-\Lambda _T^{r (X^a)}} U_{X^a,G} (X_T^a) = \mbox{} & U_{X^a,G} (x)
+ \int _0^T e^{-\Lambda _t^{r (X^a)}} {\mathcal L}_X U_{X^a,G} (X_t^a)
\, dt \nonumber \\
& + \int _0^T e^{-\Lambda _t^{r (X^a)}} U_{X^a,G}' (X_t^a) \,
dZ_t^a + M_T^a \nonumber \\
= \mbox{} & U_{X^a,G} (x) - \int _0^T e^{-\Lambda _t^{r (X^a)}}
G(X_t^a) \, dt + M_T^a , \label{Ito0}
\end{align}
where
\be
M_T^a = \int _0^T e^{-\Lambda _t^{r (X^a)}} \sigma (X_t^a) U_{X^a,G}'
(X_t^a) \, dW_t .
\ee
It\^{o}'s isometry and the continuity of $\sigma$ and
$U_{X^a,G}'$ imply that
\begin{align}
\e \left[ \left( M_{T_n^a}^a \right) ^2 \right] & = \e \left[
\int _0^\infty {\bf 1} _{\{ t \leq T_n^a \}} \left[ e^{-\Lambda
_t^{r (X^a)}} \sigma (X_t^a) U_{X^a,G}' (X_t^a) \right] ^2 dt
\right] \nonumber \\
& \leq \frac{1}{2r_0} \sup _{s \in [a,n]} \left[ \sigma (s)
U_{X^a,G}' (s) \right] ^2 \nonumber \\
& < \infty , \nonumber
\end{align}
where $r_0$ is as in Assumption~\ref{Ar}.
It follows that the stopped process $(M^a)^{T_n^a}$ is a uniformly
square integrable martingale, and therefore,
$\e \left[ (M^a)_\infty^{T_n^a} \right] \equiv \e \left[ M_{T_n^a}^a
\right] = 0$.
In view of this observation and (\ref{Ito0}), we obtain
\ben
\e \left[ e^{-\Lambda _{T_n^a}^{r (X^a)}} \right] U_{X^a,G} (n) =
U_{X^a,G} (x) - \e \left[ \int _0^{T_n^a} e^{-\Lambda _t^{r (X^a)}}
G(X_t^a) \, dt \right] . \label{Refl-1}
\een

Now, if $G$ is bounded, then (\ref{UG}) follows immediately by passing
to the limit $n \rightarrow \infty$ in (\ref{Refl-1}) using the
dominated convergence theorem, the fact that the restriction of
$U_{X^a,G}$ in $[a , \infty [$ is bounded and the assumption that
$r(x) \geq r_0 > 0$, for all $x>0$.
If $G$ is positive, we define $G^{(m)} = G \wedge m$, for $m \geq 1$,
and we note that
\ben
U_{X^a,G^{(m)}} (x) \equiv - \frac{R_{X^0,G^{(m)}}' (a)}{\phi ' (a)}
\phi (x) + R_{X^0,G^{(m)}} (x) = \e \left[ \int _0^\infty e^{-\Lambda
_t^{r (X^a)}} G^{(m)} (X_t^a) \, dt \right] \label{Refl-2}
\een
because $G^{(m)}$ is bounded.
In view of (\ref{RG}), the calculation
\be
R_{X^0,G} ' (x) = \frac{2}{C} \phi ' (x) \int_0^x \frac{G(s) \psi
(s)}{\sigma^2 (s) p_{X^0}' (s)} \, ds + \frac{2}{C} \psi ' (x) \int
_x^\infty \frac{G(s) \phi (s)} {\sigma^2(s) p_{X^0}' (s)} \, ds ,
\ee
and the monotone convergence theorem, we can see that
\be
\lim _{m \rightarrow \infty} R_{X^0,G^{(m)}} (x) = R_{X^0,G} (x) \quad
\text{and} \quad \lim _{m \rightarrow \infty} R_{X^0,G^{(m)}}' (a) =
R_{X^0,G}' (a)
\ee
because $G$ satisfies (\ref{IC-X2}), and that
\be
\lim _{m \rightarrow \infty} \e \left[ \int _0^\infty e^{-\Lambda
_t^{r (X^a)}} G^{(m)} (X_t^a) \, dt \right] = \e \left[ \int _0^\infty
e^{-\Lambda _t^{r (X^a)}} G(X_t^a) \, dt \right] .
\ee
However, these limits and (\ref{Refl-2}) imply (\ref{UG}).
The general case now follows by considering the minimal decomposition
$G = G^+ - G^-$ of the function $G$ to the difference of two positive
functions, and by linearity.

To complete the proof, we note that, if $G$ satisfies (\ref{IC-X2}),
then $U_{X^a, |G|} (x) < \infty$, which, combined with (\ref{UG}),
implies (\ref{IC-Xa}).
Also, (\ref{Refl-1}), (\ref{IC-Xa}), (\ref{UG}) and the dominated
convergence theorem imply that
$\lim _{n \rightarrow \infty} \e \left[ e^{-\Lambda _{T_n^a}^{r (X^a)}}
\right] U_{X^a,G} (n) = 0$.
However, this limit, the fact that the restriction of $\phi$ in
$[a, \infty [$ is bounded and Assumption~\ref{A3}.(c) imply
(\ref{TVC-refl}).
\mbox{}\hfill$\Box$

%=============================================================
\section{The solution to the control problem}
\label{sec:soltn}

With regard to standard theory of stochastic control (e.g., see
Chapter~VIII in Fleming and Soner~\cite{FS}), we expect that
the value function $v$ identifies with an appropriate solution $w$ to
the Hamilton-Jacobi-Bellman (HJB) equation
\ben
\max \left\{ {\mathcal L}_X w(x) + h(x) , \ w'(x) - k(x) \right\} =
0 , \label{HJB}
\een
where ${\mathcal L}_X$ is the operator defined by (\ref{L}).
We conjecture that the optimal strategy of our problem can take one of
two qualitatively different forms, depending on the problem data.
The first of these arises when it is never optimal to exert any
control effort.
In this case, we expect that the function $R_{X^0,h}$ that is
defined by (\ref{Rh-def}) in Remark~\ref{rem:RhQ} should satisfy
(\ref{HJB}).

The second one is characterised by a boundary point $a>0$ and can be
described as follows.
If the system's initial condition $x$ is less than $a$, then maximal
control should be exercised to immediately reposition the system's
state at level $a$ (i.e., cause a ``jump'' of size $a-x$ in the
positive direction at time 0).
After this possible initial jump, minimal control should be exercised
so that the system's state process is reflected at the boundary point
$a$ in the positive direction.
In view of the heuristic arguments that explain the structure of the
HJB equation (\ref{HJB}), if this strategy is indeed optimal, then
we should look for a function $w$ and a point $a>0$ such that
\ben
{\mathcal L}_X w (x) + h(x) = 0 , \quad \text{for } x \geq a ,
\label{ODE}
\een
and
\ben
w(x) = w(a) - \int _x^a k(s) \, ds , \quad \text{for } x < a .
\label{w-x<a}
\een
Now, every solution to the ODE in (\ref{ODE}) is given by
\ben
w(x) = A \phi (x) + B \psi (x) + R_{X^0,h} (x) , \label{ODE-sol}
\een
for some constants $A , B \in \R$, where the functions $\phi$ and
$\psi$ are defined as in Section~\ref{sec:ODE}.
It turns out that the arguments that we use to establish Theorem
\ref{main-thm} below, which is our main result, remain valid only
for the choice $B=0$.
For this reason, we look for a solution to the HJB
equation (\ref{HJB}) of the form
\ben
w(x) = \begin{cases} A \phi (x) + R_{X^0,h} (x) , & \text{for } x
\geq a , \\ w(a) - \int _x^a k(s) \, ds , & \text{for } x < a .
\end{cases} \label{w}
\een

To specify the parameter $A$ and the free-boundary point $a$, we
appeal to the so-called ``smooth pasting'' condition of singular
stochastic control that requires that the value function should be
$C^2$, in particular, at the free-boundary point $a$.
This requirement gives rise to the system of equations
\begin{align}
R_{X^0,h} '(a) + A \phi '(a) & = k(a) , \nonumber \\
R_{X^0,h} ''(a) + A \phi ''(a) & = k'(a) , \nonumber
\end{align}
which is equivalent to
\begin{gather}
A = \frac{k(a) - R_{X^0,h} '(a)}{\phi '(a)} = \frac{k'(a) -
R_{X^0,h} ''(a)}{\phi ''(a)} , \label{A} \\
\phi' (a) \left[ R_{X^0,h}'' (a) - k' (a) \right] - \phi ''
(a) \left[ R_{X^0,h}' (a) - k(a) \right]= 0 .
\label{wronskian=0}
\end{gather}
Taking note of (\ref{Rh-K/Q}) in Remark~\ref{rem:RhQ}, we
can see that
\ben
R_{X^0,h}' (x) - k(x) = R_{X^0,Q}' (x) , \label{Rh-K/Q'}
\een
so (\ref{wronskian=0}) is equivalent to
\ben
\phi' (a) R_{X^0,Q}'' (a) - \phi '' (a) R_{X^0,Q}' (a)
= 0 . \label{EQN}
\een
Now, Proposition~\ref{R-(X,G)'=R-(Y,H)} and (\ref{Rdrg}) imply
that
\begin{align}
R_{X^0,Q}' (x) & = R_{Y^0, \drq}(x) \nonumber \\
& = - \frac{1}{C} \phi' (x) \int _0^x \frac{(\drq) (s) \psi' (s)}
{r(s) p_{X^0}'(s)} \, ds - \frac{1}{C} \psi' (x) \int _x^\infty
\frac{(\drq) (s) \phi' (s)} {r(s) p_{X^0}'(s)} \, ds , \label{RQ'}
\end{align}
where the function $\drq$ is as in Assumption~\ref{A3}.(d).
Combining this observation with the calculation
\ben
R_{X^0,Q}'' (x) = - \frac{1}{C} \phi'' (x) \int_0^x \frac{(\drq)
(s) \psi' (s)} {r(s) p_{X^0}'(s)} \, ds - \frac{1}{C} \psi'' (x) \int
_x^\infty \frac{(\drq) (s) \phi' (s)} {r(s) p_{X^0}'(s)} \, ds ,
\label{RQ''}
\een
we can see that (\ref{EQN}) is equivalent to
\ben
\left[ \phi'' (a) \psi' (a) - \phi' (a) \psi'' (a) \right]
\int _a^\infty \frac{(\drq) (s) \phi' (s)} {r(s) p_{X^0}'(s)}
\, ds = 0 . \label{wronskian-f-phi}
\een
In view of the fact that $\phi''\psi' - \phi'\psi''$ is strictly
positive (see also (\ref{phi-psi'-Wr})), we conclude that
(\ref{wronskian=0}) is equivalent to
\ben
g(a) := \int _a^\infty \frac{(\drq) (s) \phi' (s)} {r(s) p_{X^0}'(s)}
\, ds = 0 . \label{g-defn}
\een

The following result is concerned with the considerations above
regarding the solvability of the HJB equation (\ref{HJB}).

%-------------------------------------------------------------
\begin{lem}\label{x-star-lemma}
Suppose that Assumptions~\ref{A1}--\ref{A3} are
satisfied.
The equation $g(a)=0$ has a unique solution $a>0$ if and only if
\ben
x^* > 0 \quad \text{and} \quad \lim _{x \downarrow 0} g(x) < 0
, \label{g-cond}
\een
where $x^* \geq 0$ is as in Assumption~\ref{A3}.(d).
Furthermore, the following cases hold true:
\vspace{2mm}

\noindent
(I) If (\ref{g-cond}) is not true, then the function $R_{X^0,h}$
satisfies the HJB equation (\ref{HJB}).
\vspace{2mm}

\noindent
(II) If (\ref{g-cond}) is true, $a>0$ is the unique solution to
the equation $g(a)=0$ and $A \geq 0$ is the constant given by
(\ref{A}), then the function $w$ defined by (\ref{w}) is $C^2$ and
solves the HJB equation (\ref{HJB}).
\vspace{2mm}

\noindent
In either case, the associated solution $w$ to the HJB equation
(\ref{HJB}) is bounded from below.
\label{lem:HJB} \end{lem}
{\bf Proof}.
%===========
In view of the definition (\ref{g-defn}) of $g$ and
Assumption~\ref{A3}.(d), we calculate
\ben
g'(x)= - \frac{(\drq) (x) \phi'(x)}{r(x) p_{X^0}'(x)}
\begin{cases}
\geq 0 , & \text{for } x \leq x^*, \text{ if } x^* > 0 ,
\\ < 0 & \text{for }  x > x^* .
\end{cases} \label{g'}
\een
Combining this calculation with the fact that $\lim _{x
\rightarrow \infty} g(x)=0$, we can see that the equation $g(a)
= 0$ has a unique solution $a>0$ if and only if (\ref{g-cond})
is true.
For future reference, we also note that
\ben
\text{if there exists } a>0  \text{ such that } g(a) = 0 ,
\text{ then } a < x^* \text{ and } g(x) > 0 , \text{ for
all } x>a . \label{a<x*}
\een
In this case, (\ref{g'}) and the fact that $\phi' < 0 <
\psi'$ imply that
\be
\frac{(\drq) (x) \psi'(x)}{r(x) p_{X^0}'(x)} \geq 0 ,
\quad \text{for all } x \leq a ,
\ee
which, combined with (\ref{A}), (\ref{Rh-K/Q'}), (\ref{RQ'}) and
(\ref{g-defn}), implies that
\begin{align}
A & = - \frac{R_{X^0,Q} '(a)}{\phi' (a)} = \frac{1}{C} \int _0^a
\frac{(\drq) (s) \psi' (s)} {r(s) p_{X^0}'(s)} \, ds
\geq 0 . \nonumber
\end{align}
On the other hand,
\ben
\text{if the equation } g(a) = 0 \text{ has no solution }
a>0 , \text{ then } g(x) > 0 , \text{ for all } x>0 .
\label{g(a)-nosol}
\een

With regard to Case~I we will prove that $R_{X^0, h}$
satisfies the HJB equation (\ref{HJB}) if we show that
$R_{X^0, h}' (x) - k(x) \leq 0$, for all $x>0$, which
is equivalent to showing that
\ben
R_{Y^0,\drq} (x) \leq 0 , \quad \text{for all } x>0 ,
\label{I-HJB}
\een
thanks to (\ref{Rh-K/Q'}) and (\ref{RQ'}).
To this end, we use (\ref{RQ'}) and (\ref{RQ''}) to
calculate
\begin{align}
\frac{d}{dx} \left( \frac{R_{Y^0,\drq}(x)}{- \phi' (x)} \right)
= - \frac{\phi'' (x) \psi' (x) - \phi' (x) \psi'' (x)}{C [\phi'
(x)]^2} g(x) . \label{RY/phi'}
\end{align}
The right-hand side of this identity is strictly negative
for all $x>0$, thanks to the strict positivity of the function
$\phi''\psi' - \phi'\psi''$ (see (\ref{phi-psi'-Wr}))
and (\ref{g(a)-nosol}).
Also, in view of (\ref{g(a)-nosol}) and the fact that $\phi' <
0 < \psi'$, we can see that
\begin{align}
\lim _{x \downarrow 0} \frac{R_{Y^0, \drq} (x)} {-\phi '(x)}
& = \frac{1}{C} \lim _{x \downarrow 0} \left[ \int _0^x  \frac{(\drq)
(s) \psi'(s)} {r(s) p_{X^0} '(s)} ds + \frac{\psi '(x)}
{\phi' (x)} g(x) \right] \nonumber \\
& \leq 0 . \nonumber
\end{align}
However, these observations imply (\ref{I-HJB}).

By construction, Case~II will follow if we show that
\begin{align}
\half \sigma ^2 (x) k'(x) + b(x) k(x) - r(x) \left[ w(a)
- \int _x^a k(s) \, ds \right]  + h(x) & \leq 0 , \quad
\text{for all } x \leq a , \label{HJB-in1} \\
w'(x) - k(x) & \leq 0 , \quad \text{for all } x > a .
\label{HJB-in2}
\end{align}
In view of the definitions (\ref{K}) and (\ref{Q}) of the
function $K$ and $Q$, respectively, we can see that
(\ref{HJB-in1}) is equivalent to
\be
Q(x) + r(x) K(a) - r(x) w(a) \leq 0, \quad \text{for all }
x \leq a .
\ee
Also, the facts that $w$ is $C^2$ and satisfies (\ref{ODE})
imply that
\begin{align}
r(a) w(a) & = \half \sigma ^2 (a) k'(a) + b(a) k(a) + h(a)
\nonumber \\
& = Q(a) + r(a) K(a) . \nonumber
\end{align}
These observations and the strict positivity of $r$ imply
that (\ref{HJB-in1}) is equivalent to
\be
\frac{Q(x)}{r(x)} \leq \frac{Q(a)}{r(a)} , \quad
\text{for all } x \leq a  .
\ee
However, this inequality follows immediately from the fact
that (\ref{a<x*}) and Assumption~\ref{A3}.(d) imply that
$\left( Q/r \right) ' (x) \geq 0$, for all $x \leq a$.

Substituting the first expression in (\ref{A}) for
$A$ in (\ref{w}), and using (\ref{Rh-K/Q'}) and (\ref{RQ'}),
we can see that (\ref{HJB-in2}) is equivalent to
\be
\frac{R_{Y^0,\drq} (a)}{-\phi '(a)} \geq \frac{R_{Y^0,\drq}
(x)} {-\phi'(x)}, \quad \text{for all } x > a .
\ee
However, this inequality follows immediately once we combine
(\ref{RY/phi'}) with the strict positivity of $\phi''\psi'
- \phi'\psi''$ and (\ref{a<x*}) to obtain
\be
\frac{d}{dx} \left( \frac{R_{Y^0,\drq}(x)}{- \phi' (x)} \right)
< 0 , \quad \text{for all } x>a .
\ee

Finally, we note that the assumption that $h/r$ is bounded from below
(see Assumption~\ref{A3}.(a)) and (\ref{Rh-def}) in
Remark~\ref{rem:RhQ} imply that $R_{X^0,h}$ is bounded from
below.
However, this observation and the structure of the solution $w$ to the
HJB equation (\ref{HJB}) associated with either of the two cases
considered imply that $w$ is bounded from below.
\mbox{}\hfill$\Box$
\vspace{5mm}

We can now prove the main result of the paper.

%-------------------------------------------------------------
\begin{thm}
Suppose that Assumptions~\ref{A1}--\ref{A3} hold.
The value function $v$ of our control problem identifies with the
solution $w$ to the HJB equation (\ref{HJB}) derived in
Lemma~\ref{lem:HJB}.
In particular, the following two cases hold true:
\vspace{2mm}

\noindent
(I) If the problem data are such that (\ref{g-cond}) is false, then it
is optimal to exert no control effort at all times.
\vspace{2mm}

\noindent
(II) If the problem data are such that (\ref{g-cond}) is true, then
the optimal intervention strategy involves a jump of size $(a-x)^+$ at
time 0 and then reflects the state process $X$ in the boundary point
$a>0$ in the positive direction.
\label{main-thm} \end{thm}
{\bf Proof.}
%===========
Fix any initial condition $x>0$ and any admissible intervention
strategy $Z \in {\mathcal A}$.
Using It\^{o}'s formula and the fact that $\Delta X_t = \Delta Z_t$, we
calculate
\begin{align}
e^{-\Lambda _T^{r (X)}} w(X_{T+}) = \mbox{} & w(x)  + \int _0^T
e^{-\Lambda _t^{r (X)}} {\mathcal L}_X w(X_t) \, dt + \int _0^T
e^{-\Lambda _t^{r (X)}} w'(X_t) \, dZ_t  \nonumber \\
& + \sum _{0 \leq t \leq T} e^{-\Lambda _t^{r (X)}} \left[ w(X_{t+}) -
w(X_t) - w'(X_t) \Delta X_t \right] + M_T \nonumber\\
= \mbox{} & w(x)  + \int _0^T e^{-\Lambda _t^{r (X)}} {\mathcal L}
_X w(X_t) \, dt + \int _0^T e^{-\Lambda _t^{r (X)}}w'(X_t) \, dZ_t^c
\nonumber \\
& + \sum _{0 \leq t \leq T} e^{-\Lambda _t^{r (X)}} \left[ w(X_t + \Delta
Z_t) - w(X_t) \right] + M_T , \nonumber
\end{align}
where the operator ${\mathcal L}_X$ is defined by (\ref{L}), and
\be
M_T = \int _0^T e^{-\Lambda _t^{r (X)}}\sigma (X_t) w'(X_t) \, dW_t .
\ee
In view of this calculation, (\ref{Z-integral}) and the fact that $w$
satisfies the HJB equation (\ref{HJB}), we obtain
\begin{align}
\int _0^T e^{-\Lambda_t^{r (X)}} h(X_t) \, & dt - \int _0^T
e^{-\Lambda _t^{r (X)}} k(X_t) \circ dZ_t \nonumber \\
= \mbox{} & w(x) - e^{-\Lambda _T^{r (X)}} w(X_{T+}) + \int _0^T
e^{-\Lambda _t^{r (X)}} \left[ {\mathcal L}_X w(X_t) + h(X_t) \right]
dt \nonumber \\
& + \int _0^T e^{-\Lambda _t^{r (X)}} \left[ w'(X_t) - k(X_t) \right]
dZ_t^c \nonumber \\
& + \sum _{0 \leq t \leq T} e^{-\Lambda _t^{r (X)}} \int _0^{\Delta
Z_t} \left[ w'(X_t+s) - k(X_t+s) \right] ds + M_T \nonumber \\
\leq \mbox{} & w(x) - e^{-\Lambda _T^{r (X)}} w(X_{T+}) + M_T
. \label{Ito1}
\end{align}

To proceed further, let $( \tau_n)$ be a localising sequence for the
stochastic integral $M$.
Taking expectations in (\ref{Ito1}), we obtain
\begin{align}
\e \biggl[ \int _0^{\tau _n \wedge T} e^{-\Lambda _t^{r (X)}} h(X_t)
\, dt - \int _0^{\tau _n \wedge T} e^{-\Lambda _t^{r (X)}} k(X_t)
\circ dZ_t \biggr] \leq w(x) + \e \left[ e^{-\Lambda_{\tau _n \wedge
T}^{r (X)}} w^- \left( X_{( \tau _n \wedge T) +} \right) \right] ,
\label{Ito3}
\end{align}
where $w^- = - (w \wedge 0)$.
Now, the assumption that $h/r$ is bounded from below, the fact that
the process $\Theta$ defined by $\Theta_t = - \exp \left( - \Lambda
_t^{r(X)} \right)$, for $t \geq 0$, is increasing, and Fatou's lemma
imply that
\begin{align}
\e \left[ \int _0^T e^{-\Lambda _t^{r (X)}} h(X_t) \, dt \right]
& = \e \left[ \int _0^T \frac{h(X_t)}{r(X_t)} \, d\Theta_t \right]
\nonumber \\
& \leq \liminf _{n \rightarrow \infty} \e \left[ \int _0^T {\bf 1}
_{\{ t \leq \tau_n \}} \frac{h(X_t)}{r(X_t)} \, d\Theta_t \right]
\nonumber \\
& = \liminf _{n \rightarrow \infty} \e \left[ \int _0^{\tau _n \wedge T}
e^{-\Lambda _t^{r (X)}} h(X_t) \, dt \right] ,
\end{align}
while the monotone convergence theorem implies that
\be
\e \left[ \int _0^T e^{-\Lambda _t^{r (X)}} k(X_t) \circ dZ_t
\right] = \lim _{n \rightarrow \infty} \e \left[ \int _0^{\tau _n
\wedge T} e^{-\Lambda _t^{r (X)}} k(X_t) \circ dZ_t \right] ,
\ee
Also, the fact that $w^-$ is bounded (see Lemma~\ref{lem:HJB}) and the
dominated convergence theorem imply that
\be
\lim _{n \rightarrow \infty} \e \left[ e^{-\Lambda _{\tau _n \wedge
T}^{r (X)}} w^- (X_{( \tau _n \wedge T) +}) \right] = \e \left[
e^{-\Lambda _T^{r (X)}} w^- (X_{T+}) \right] .
\ee
In view of these observations, we can pass to the limit $n \rightarrow
\infty$ in (\ref{Ito3}) to obtain
\be
\e \left[ \int _0^T e^{-\Lambda_t^{r (X)}} h(X_t) \, dt - \int
_0^T e^{-\Lambda _t^{r (X)}} k(X_t) \circ dZ_t \right] \leq w(x)
+ \e \left[ e^{-\Lambda _T^{r (X)}} w^- (X_{T+}) \right] .
\ee
Combining this inequality with the limit
\be
\lim _{T \rightarrow \infty} \e \left[ e^{-\Lambda _T^{r (X)}}
w^- (X_{T+}) \right] = 0 ,
\ee
which follows from Assumption~\ref{A3}.(c) and the fact that $w^-$ is
bounded, we can see that
\ben
J_x (Z) \equiv \limsup _{T \rightarrow \infty} \e \left[ \int _0^T
e^{-\Lambda _t^{r (X)}} h(X_t) \, dt - \int _0^T e^{-\Lambda
_t^{r (X)}} k(X_t) \circ dZ_t \right]\leq w(x). \label{VT1}
\een

If (\ref{g-cond}) is false, then the control strategy $Z^0 \equiv 0$
has payoff
\be
J_x (Z^0) = R_{X^0,h} (x) = w(x) .
\ee
The first equality here follows from the definition (\ref{Rh-def})
of $R_{X^0,h}$ and the dominated convergence theorem (see
(\ref{IC1}) in Assumption \ref{A3}.(e)), while the second one is
just Case~I of Lemma~\ref{lem:HJB}.
However, these identities and (\ref{VT1}) establish Case~I of the
theorem.

Now, let us assume that (\ref{g-cond}) is true.
Let $Z^a$ be the c\`{a}gl\`{a}d process that has a jump of size
$(a-x)^+$ at time 0 and then reflects the state process in the
boundary point $a>0$ in the positive direction (see the discussion
preceding Lemma~\ref{refl-lemma}).
Also, let $X^a$ be the associated solution to the SDE (\ref{SDE-X})
and let $T_n^a$ be the first hitting time of $\{ n \}$, which is given
by (\ref{Tan}) in Lemma~\ref{refl-lemma}.
In this case, we can check that (\ref{Ito1}) holds with equality, so
\ben
\int _0^T e^{-\Lambda_t^{r (X^a)}} h(X_t^a) \, dt - \int _0^T
e^{-\Lambda _t^{r (X^a)}} k(X_t^a) \circ dZ_t^a = w(x) - e^{-\Lambda
_T^{r (X^a)}} w(X_T^a) + M_T^a , \label{Ito2}
\een
for all $T>0$, where
\be
M_T^a = \int _0^T e^{-\Lambda _t^{r (X^a)}} \sigma (X_t^a) w'(X_t^a)
\, dW_t .
\ee
Using It\^{o}'s isometry and the fact that $\sigma$ and $w'$ are
continuous, we can see that, for $n \geq a \vee x$,
\begin{align}
\e \left[ \left( M_{T_n^a}^a \right) ^2 \right] & = \e \left[
\int _0^\infty {\bf 1} _{\{ t \leq T_n^a \}} \left[ e^{-\Lambda
_t^{r (X^a)}} \sigma (X_t^a) w'(X_t^a) \right] ^2 dt \right]
\nonumber \\
& \leq \frac{1}{2r_0} \sup _{s \in [a,n]} \left[ \sigma (s)
w'(s) \right] ^2 \nonumber \\
& < \infty , \nonumber
\end{align}
where $r_0$ is as in Assumption~\ref{Ar}.
It follows that the stopped process $\left( M^a \right) ^{T_n^a}$ is a
uniformly square integrable martingale, so $\e \left[ \left( M^a
\right) _\infty^{T_n^a} \right] \equiv \e \left[ M_{T_n^a}^a \right] =
0$.
This observation and (\ref{Ito2}) imply that
\ben
\e \left[ \int _0^{T_n^a} e^{-\Lambda _t^{r (X^a)}} h(X_t^a) \, dt -
\int _0^{T_n^a} e^{-\Lambda _t^{r (X^a)}} k(X_t^a) \circ dZ_t^a
\right] = w(x) - \e \left[ e^{-\Lambda_{T_n^a}^{r (X^a)}} \right]
w(n) , \label{Ito13}
\een
for all $n > a \vee x$.

Combining the fact that $w$ is the sum of $R_{X^0,h}$ and a bounded
function with Assumption~\ref{A3}.(c) and (\ref{TVC-refl}) in Lemma
\ref{refl-lemma},
we can see that
\be
\lim _{n \rightarrow \infty} \e \left[ e^{-\Lambda _{T_n^a}
^{r (X^a)}} \right] w(n) = 0 .
\ee
Also, Lemma~\ref{refl-lemma} and the dominated convergence theorem
imply that
\begin{align}
\lim _{n \rightarrow \infty} \e \left[ \int _0^{T_n^a} e^{-\Lambda
_t^{r (X^a)}} h(X_t^a) \, dt \right] = \e \left[ \int _0^\infty
e^{-\Lambda _t^{r (X^a)}} h(X_t^a) \, dt \right] \in \R , \label{VT21}
\end{align}
while the monotone convergence theorem implies that
\be
\lim _{n \rightarrow \infty} \e \left[ \int _0^{T_n^a}
e^{-\Lambda _t^{r (X^a)}} k(X_t^a) \circ dZ_t^a \right] = \e \left[
\int _0^\infty e^{-\Lambda _t^{r (X^a)}} k(X_t^a) \circ dZ_t^a
\right] .
\ee
These limits and (\ref{Ito13}) imply that
\be
\e \left[ \int _0^\infty e^{-\Lambda _t^{r (X^a)}} h(X_t^a) \, dt
\right] - \e \left[ \int _0^\infty e^{-\Lambda _t^{r (X^a)}}
k(X_t^a) \circ dZ_t^a \right] = w(x) .
\ee
If we combine this conclusion with (\ref{VT21}), then we can see that
\be
\e \left[ \int _0^\infty e^{-\Lambda _t^{r (X^a)}} k(X_t^a) \circ
dZ_t^a \right] < \infty ,
\ee
so $Z^a$ is admissible.
Furthermore, if we combine it with (\ref{VT1}), then we can see that
$v(x) = w(x)$ and that $Z^a$ is optimal.
\mbox{}\hfill$\Box$

%=============================================================
\section{Special cases}
\label{sec:examples}

We now consider special cases that arise when the running payoff
function $h$ is a power utility function and the running cost function
$k$ as well as the discounting factor $r$ are constant.
In particular, we assume that
\be
h(x) = \lambda x^\nu , \quad k(x) = \kappa \quad \text{and} \quad r(x)
= r_1 ,
\ee
for some constants $\kappa , \lambda , r_1 > 0$ and $\nu \in \,
]0,1[$.
Also, we assume that the uncontrolled system's state dynamics are
modelled by a geometric Brownian motion (Section~\ref{sec:GBMex}) or
by a mean-reverting square-root process such as the one in the
Cox-Ingersoll-Ross interest rate model (Section~\ref{sec:CIRex}).

%=============================================================
\subsection{Geometric Brownian motion}
\label{sec:GBMex}

Suppose that $X^0$ is a geometric Brownian motion, so that
\be
dX_t^0 = b X_t^0 \, dt + \sigma X_t^0 \, dW_t , \quad X_0^0 = x > 0 ,
\ee
for some constants $b$ and $\sigma \neq 0$, and assume that $r_1 > b$.
In this case, it is a standard exercise to verify that, if we choose
$c=1$, then
\be
\phi (x) = x^m , \quad \psi (x) = x^n \quad \text{and} \quad p_{X^0}'
(x) = x^{n+m-1} ,
\ee
where the constants $m<0<n$ are the solution to the quadratic equation
\be
\half \sigma ^2 l^2 + \left( b - \half \sigma ^2 \right) l - r_1 = 0 .
\ee
Also, it is well-known that
\ben
r_1 > b \quad \Leftrightarrow \quad n>1 . \label{GBM-rb}
\een

In this context, Assumptions~\ref{A1}, \ref{A2}, \ref{Ar}, \ref{AY}
and~\ref{A3}.(a)--\ref{A3}.(c) are plainly satisfied.
Also, in the presence of (\ref{GBM-rb}), we can calculate
\be
\e \left[ \int _0^\infty e^{-r_1 t} X_t \, dt \right] =
\frac{x}{r_1-b} .
\ee
However, using this observation and the facts that
\be
0 \leq h(x) \leq \lambda (1+x) \quad \text{and} \quad {\mathcal L}_X
K(x) = - (r_1-b) \kappa x ,
\ee
we can verify that Assumption~\ref{A3}.(e) is satisfied.
Furthermore, we can use the calculation
\be
\drq (x) = Q'(x) = \lambda \nu x^{-(1-\nu)} - \kappa (r_1-b)
\ee
to verify that $x^* > 0$ and that Assumption~\ref{A3}.(d) is satisfied
as well.

In view of the fact that $n>1$, we can check that the function $g$
defined by (\ref{g-defn}) admits the expression
\be
g(x) = \frac{m \lambda \nu}{n - \nu} x^{\nu - n} - \frac{m \kappa
(r_1-b)} {n-1} x^{1-n} , \quad \text{for } x > 0 .
\ee
It follows that Case~II of Theorem~\ref{main-thm} is always true
and that the unique solution $a>0$ to the equation $g(a)
= 0$ that characterises the optimal strategy is given by
\be
a = \left[ \frac{\lambda \nu (n-1)} {\kappa (r_1-b) (n - \nu )}
\right] ^{1 / (1-\nu )} .
\ee

%=============================================================
\subsection{Mean-reverting square-root process}
\label{sec:CIRex}

Suppose that $X^0$ is a mean-reverting square-root process, so that
\be
dX_t^0 = \alpha (\theta -X_t^0) \, dt + \sigma \sqrt{X_t^0}  \, dW_t ,
\quad X_0^0 = x > 0 ,
\ee
for some constants $\alpha , \theta , \sigma > 0$, and assume that
\ben
\alpha \theta - \half \sigma ^2 > 0 , \label{CIR-cond}
\een
which is a necessary and sufficient condition for $X^0$ to be
non-explosive.
In this case, we can check that the associated diffusion $Y^0$
satisfies the SDE
\be
dY_t^0 = \alpha \left( \left[ \theta + \frac{\sigma ^2}{2\alpha}
\right] - Y_t^0 \right) dt + \sigma \sqrt{Y_t^0}  \, dW_t , \quad
Y_0^0 = x > 0 ,
\ee
and that $Y^0$ is non-explosive because $\alpha \theta > 0$, which is
plainly true, is the condition corresponding to (\ref{CIR-cond}).
We can therefore see that Assumptions~\ref{A1}, \ref{A2}, \ref{Ar},
\ref{AY} and \ref{A3}.(a)--\ref{A3}.(c) are all satisfied.

Now, we calculate
\begin{align}
\e \left[ \int _0^\infty e^{-r_1 t} X_t \, dt \right] & = \int
_0^\infty e^{-r_1 t} \left[ \theta + (x-\theta) e^{-\alpha t} \right]
dt \nonumber \\
& = \frac{\alpha \theta + r_1 x}{r_1 (\alpha + r_1)} . \nonumber
\end{align}
Combining these identities with the inequalities
\be
0 \leq h(x) \leq \lambda (1+x) \quad \text{and} \quad \left|
{\mathcal L}_X K(x) \right| = \left| \alpha \theta \kappa - \kappa
(\alpha + r_1) x \right| \leq \alpha \theta \kappa + \kappa
(\alpha + r_1) x ,
\ee
we can verify that Assumption~\ref{A3}.(e) is satisfied.
Also, in view of the calculation
\be
\drq (x) = Q'(x) = \lambda \nu x^{-(1-\nu)} - (\alpha + \kappa r_1) ,
\ee
we can see that $x^* > 0$ and that Assumption~\ref{A3}.(d) holds.

Making the transformations $y = 2 \alpha x / \sigma ^2$ and $\hat{w}
(y) = w (x)$ in the ODE ${\mathcal L}_X w(x) =0$, we obtain
\be
\hat{w}'' (y) + \left( \frac{2\alpha \theta}{\sigma ^2} - y \right)
\hat{w}' (y)- \frac{r_1}{\alpha} \hat{w} (y) = 0 ,
\ee
which is Kummer's equation.
With reference to Abramowitz and Stegun~\cite[Chapter~13]{AS} or
Magnus, Oberhettinger and Soni~\cite[Chapter VI]{MOS}, it follows
that, if we choose $c=1$, then
\begin{gather}
\phi (x) = \frac{U \left( \frac{r_1}{\alpha} , \frac{2\alpha \theta}
{\sigma ^2} ; \frac{2\alpha}{\sigma ^2} x \right)} {U \left(
\frac{r_1}{\alpha} , \frac{2\alpha \theta}{\sigma ^2} ;
\frac{2\alpha}{\sigma ^2} \right)} , \quad
\psi(x) = \frac{\mbox{}_1F_1 \left( \frac{r_1}{\alpha} , \frac{2\alpha
\theta} {\sigma ^2} ; \frac{2\alpha}{\sigma ^2} x \right)}
{\mbox{}_1F_1 \left( \frac{r_1}{\alpha} , \frac{2\alpha \theta} {\sigma
^2} ; \frac{2\alpha}{\sigma ^2} \right)} \nonumber \\
\intertext{and}
p_{X^0}' (x) = x^{-2 \alpha \theta / \sigma ^2} e^{2\alpha (x-1) /
\sigma ^2} . \nonumber
\end{gather}
Note that $\mbox{}_1F_1$ is the well-known confluent hypergeometric
function.

Now, in view of the differentiation formula
\be
U'(a,b;x) = - a U(a+1,b+1;x)
\ee
(see Abramowitz and Stegun~\cite[13.4.21]{AS}), we can see that the
function $g$ defined by (\ref{g-defn}) is given by
\begin{align}
g(x) = \mbox{} - \frac{2r_1 e^{2\alpha / \sigma ^2}} {U \left(
\frac{r_1}{\alpha} , \frac{2\alpha \theta}{\sigma ^2} ;
\frac{2\alpha}{\sigma ^2} \right)} \biggl[ & \lambda \nu \int
_x^\infty s^{\frac{2 \alpha \theta}{\sigma ^2} + \nu -1}
e^{-\frac{2\alpha}{\sigma ^2} s} U \left( \frac{r_1}{\alpha} + 1 ,
\frac{2\alpha \theta} {\sigma ^2} + 1 ; \frac{2\alpha}{\sigma ^2}
s \right) ds \nonumber \\
& - (\alpha + \kappa r_1) \int _x^\infty s^{\frac{2 \alpha
\theta}{\sigma ^2}} e^{-\frac{2\alpha}{\sigma ^2} s} U \left(
\frac{r_1}{\alpha} + 1 , \frac{2\alpha \theta} {\sigma ^2} + 1 ;
\frac{2\alpha}{\sigma ^2} s \right) ds \biggr] . \nonumber
\end{align}
This expression is complex, so we have to resort to numerical
techniques to determine the solution of the equation $g(a) = 0$.

%=============================================================

\end{document}